\documentclass[12pt,thmsa]{article}
\usepackage{amssymb}
\begin{document}

\author{David Carf\`{i}}
\title{Schwartz families in tempered distribution spaces}
\date{}
\maketitle

\begin{abstract}
In this paper we define Schwartz families in tempered distribution spaces
and prove many their properties. Schwartz families are
the analogous of infinite dimensional matrices of separable Hilbert spaces,
but for the Schwartz test function spaces, having elements (functions)
realizable as vectors indexed by real Euclidean spaces (ordered families of
scalars indexed by real Euclidean spaces). In the paper, indeed, one of the
consequences of the principal result (the characterization of summability for Schwartz
families) is that the space of linear continuous operators among Schwartz
test function spaces is linearly isomorphic with the space of Schwartz
families. It should be noticed that this theorem is possible because of the
very good properties of Schwartz test function spaces and because of the particular
structures of the Schwartz families viewed as generalized matrices; in fact,
any family of tempered distribution, regarded as generalized matrix, has one
index belonging to a Euclidean space and one belonging to a test function
space, so that any Schwartz family is a matrix of type $(\Bbb{R}^{m},\Bbb{R%
}^{n})$ in the sense of distributions. Another motivation for the
introduction and study of these families is that these are the families
which are summable with respect to every tempered system of coefficients, in
the sense of superpositions. The Schwartz families we present in this
paper are one possible rigorous and simply manageable mathematical model for
the infinite matrices used frequently in Quantum Mechanics.
\end{abstract}

\bigskip

\section{\textbf{Families of distributions}}

\bigskip

Let $n$ be a positive natural number; by $\mathcal{S}_{n}^{\prime }$ we
denote the space of complex (or real) tempered distributions defined on the $%
n$-dimensional real Euclidean space. Let $I$ be a non-empty set, we shall
denote by $(\mathcal{S}_{n}^{\prime })^{I}$ the space of all the families in
the space of tempered distributions $\mathcal{S}_{n}^{\prime }$ indexed by
the set $I$, i.e., the set of all the surjective maps from the set $I$ onto
a subset of the space $\mathcal{S}_{n}^{\prime }$. Moreover, as usual, if $v$
is any family in the space of tempered distributions $\mathcal{S}%
_{n}^{\prime }$ and indexed by the set $I$, for each index $p\in I$, the
distribution $v(p)$ (corresponding to the index $p$ in the map $v$) is
denoted by $v_{p}$, and the family $v$ itself is also denoted by the
expressive notation $(v_{p})_{p\in I}$.

\bigskip

\textbf{Dangerous bend.} Note that it is not correct and strongly misleading
to consider the family $v$ coincident with its trace (trajectory, underlying
set) $v(I)$, i.e. its image in $\mathcal{S}_{n}^{\prime }$: the family $v$
is a function, an ordered set, and not a simple set of objects.

\bigskip

\textbf{Linear operations for families.} The set $(\mathcal{S}_{n}^{\prime
})^{I}$ of all families in the space $\mathcal{S}_{n}^{\prime }$ indexed by
a non-empty set $I$ is a vector space with respect to the following standard
binary operations:

\begin{itemize}
\item  \emph{the componentwise addition} 
\[
+:(\mathcal{S}_{n}^{\prime })^{I}\times (\mathcal{S}_{n}^{\prime
})^{I}\rightarrow (\mathcal{S}_{n}^{\prime })^{I}, 
\]
defined pointwise by 
\[
v+w:=(v_{p}+w_{p})_{p\in I}, 
\]
for any two families $v,w$;

\item  \emph{the componentwise multiplication by scalars} 
\[
\cdot \;:\Bbb{K}\times (\mathcal{S}_{n}^{\prime })^{I}\rightarrow (\mathcal{S%
}_{n}^{\prime })^{I} 
\]
defined pointwise by 
\[
av:=(av_{p})_{p\in I}, 
\]
for any family $v$ and any scalar $a$.
\end{itemize}

\bigskip

In other words, the sum family $v+w$ is defined by 
\[
\left( v+w\right) _{p}=v_{p}+w_{p}, 
\]
for every index $p$ in $I$, and the product family $av$ is defined by 
\[
\left( av\right) _{p}=av_{p}, 
\]
for every $p$ in $I$.

\bigskip

\textbf{Interpretations.} We are already familiar with ordered systems of
vectors in linear spaces: the ordered bases, which permit to define the
systems of components; the Hilbert bases of a Hilbert space, which permit to
define the abstract Fourier transforms; the Schauder bases of topological
vector spaces; and so on. But there is another natural way to think at a
(ordered) family of vectors: it can be viewed as a matrix, the unique matrix
having the ordered family as ordered systems of its rows (or columns). This
last interpretation is the useful way to think at families of distributions,
when such families should be summed or multiplied by scalars.

\bigskip

\section{\textbf{Action of a family on test functions}}

\bigskip

The basic important consideration for our purposes is the observation that a
family of tempered distributions can act naturally on test functions, as the
following definition specifies. This capacity of action is not surprising,
if we think to the families as matrices and to the test functions as vectors.

\bigskip

\textbf{Definition (image of a test function by a family of distributions).} 
\emph{Let }$v$\emph{\ be a family of tempered distributions in the space }$%
\mathcal{S}_{n}^{\prime }$\emph{\ indexed by a non-empty set }$I$\emph{\ and
let }$\phi \in \mathcal{S}_{n}$\emph{\ be any test function. The mapping} 
\[
v(\phi ):I\rightarrow \Bbb{K} 
\]
\emph{defined by} 
\[
v(\phi )(p):=v_{p}(\phi ), 
\]
\emph{for each index }$p\in I$\emph{, is called the \textbf{image of the
test function} }$\phi $\emph{\ \textbf{under the family of tempered
distributions} }$v.$

\bigskip

So, in the conditions of the above definition, thefunction $v(\phi )$
resulting from the action of the family $v$ upon the test function $\phi $
is a function belonging to the function space $\mathcal{F}(I,\Bbb{K})$
(following Bourbaki, this is the set of all scalar function on $I$).

\bigskip

Thus, with any family $v$ belonging to the vector space $(\mathcal{S}%
_{n}^{\prime })^{I}$ we can associate a mapping from the space of test
functions $\mathcal{S}_{n}$ into the function space $\mathcal{F}(I,\Bbb{K})$.

\bigskip

\textbf{Matrix interpretation.} Note that if we think to a family of
distribution $v$ belonging to the space $\mathcal{F}(I,\mathcal{S}%
_{n}^{\prime })$ as a generalized matrix with $I$ rows and $\Bbb{R}^{n}$
columns (actually the second index is belonging to the test function space $%
\mathcal{S}_{n}$, but we can think, in the sense of distributions, that it
is belonging to $\Bbb{R}^{n}$) and if we think to any test function $f$ in $%
\mathcal{S}_{n}$ as to a scalar vector with $\Bbb{R}^{n}$ components, it is
natural to understand that the result $v(f)$ of the action of the family $v$
on the function $f$ is a scalar vector with $I$ components. We so have build
up a multiplication 
\[
\mathcal{F}(I,\mathcal{S}_{n}^{\prime })\times \mathcal{S}(\Bbb{R}^{n},\Bbb{K%
})\rightarrow \mathcal{F}(I,\Bbb{K}) 
\]
analogous to that defined for finite matrices.

\bigskip

\textbf{Remark.} Equivalently, for every test function $\phi $, we have a
``projection'' $\pi _{\phi }$ sending any family of $(\mathcal{S}%
_{n}^{\prime })^{I}$ into a scalar family of the product $(\Bbb{K})^{I}$: 
\[
\pi _{\phi }(v)=(v_{p}(\phi ))_{p\in I}, 
\]
for every family $v$ of the space $(\mathcal{S}_{n}^{\prime })^{I}$.

\bigskip

\section{$^{\mathcal{S}}$\textbf{Families}}

\bigskip

In the Theory of Superpositions on the space of tempered distributions $%
\mathcal{S}_{n}^{\prime }$ the below class of $^{\mathcal{S}}$families plays
a basic role.

\bigskip

\textbf{Definition}\emph{\ }\textbf{(family of tempered distributions of
class\ }$\mathcal{S}$\textbf{).}\emph{\ Let }$v$\emph{\ be a family of
distributions in the space }$\mathcal{S}_{n}^{\prime }$\emph{\ indexed by
the Euclidean space }$\Bbb{R}^{m}$\emph{. The family }$v$\emph{\ is called a 
\textbf{Schwartz family} or\ \textbf{family of class} }$\mathcal{S}$\emph{\
or even }$^{\mathcal{S}}$\emph{\textbf{family} if, for each test function }$%
\phi \in \mathcal{S}_{n}$\emph{, the image of the test function }$\phi $%
\emph{\ by the family }$v$\emph{\ - that is the function }$v(\phi ):\Bbb{R}%
^{m}\rightarrow \Bbb{K}$\emph{\ defined by} 
\[
v(\phi )(p):=v_{p}(\phi ), 
\]
\emph{for each index }$p\in \Bbb{R}^{m}$\emph{\ - belongs to the space of
test functions }$\mathcal{S}_{m}.$\emph{\ We shall denote the set of all }$^{%
\mathcal{S}}$\emph{families by }$\mathcal{S}(\Bbb{R}^{m},\mathcal{S}%
_{n}^{\prime })$\emph{.}

\bigskip

\textbf{Example (the Dirac family in }$\mathcal{S}_{n}^{\prime }$\textbf{).}%
\emph{\ The \textbf{Dirac family in} }$\mathcal{S}_{n}^{\prime }$\emph{,
i.e., the family }$\delta :=(\delta _{x})_{x\in \Bbb{R}^{n}}$\emph{, where }$%
\delta _{x}$\emph{\ is the Dirac (tempered) distribution centered at the
point }$x$\emph{\ of }$\Bbb{R}^{n}$\emph{, is a Schwartz family.}

\emph{\bigskip }

\emph{Proof.} Indeed, for each test function $\phi \in \mathcal{S}_{n}$ and
for each index (point) $x$ in $\Bbb{R}^{n}$, we have 
\[
\delta (\phi )(x)=\delta _{x}(\phi )=\phi (x), 
\]
and hence the image of any test function $\phi $ by the Dirac family is the
function $\delta (\phi )=\phi $. So the image of the test function $\phi $
under the family $\delta $ is the function $\phi $ itself, which (in
particular) lies in the Schwartz space $\mathcal{S}_{n}$. $\blacksquare $

\bigskip

It is clear that the space of $^{\mathcal{S}}$families in $\mathcal{S}
_{n}^{\prime }$, indexed by some Euclidean space $I$, is a subspace of the
vector space $(\mathcal{S}_{n}^{\prime })^{I}$ of all families in $\mathcal{S%
}_{n}^{\prime }$ indexed by the same index set $I$.

\bigskip

\section{$^{\mathcal{S}}$\textbf{Family generated by an operator}}

\bigskip

In this section we introduce a wide class of $^{\mathcal{S}}$families. We
will see later that this class is indeed the entire class of Schwartz
families. We recall that by $\sigma (\mathcal{S}_{n})$ we denote the weak
topology $\sigma (\mathcal{S}_{n},\mathcal{S}_{n}^{\prime })$.

\bigskip

\textbf{Theorem (on the }$^{\mathcal{S}}$\textbf{family generated by a
linear and continuous operator).}\emph{\ Let }$A:\mathcal{S}_{n}\rightarrow 
\mathcal{S}_{m}$\emph{\ be a linear and continuous operator with respect to
the natural topologies of }$\mathcal{S}_{n}$\emph{\ and }$\mathcal{S}_{m}$%
\emph{\ (or equivalently, continuous with respect to the weak topologies }$%
\sigma (\mathcal{S}_{n})$\emph{\ and }$\sigma (\mathcal{S}_{m})$\emph{) and
let }$\delta $\emph{\ be the Dirac family in }$\mathcal{S}_{m}^{\prime }$%
\emph{. Then, the family of functionals} 
\[
A^{\vee }:=(\delta _{p}\circ A)_{p\in \Bbb{R}^{m}} 
\]
\emph{is a family of distribution and it is an }$^{\mathcal{S}}$\emph{family.%
}

\emph{\bigskip }

\emph{Proof.} Let $A:\mathcal{S}_{n}\rightarrow \mathcal{S}_{m}$ be a linear
and continuous operator with respect to the natural topologies of $\mathcal{S%
}_{n}$ and $\mathcal{S}_{m}$ (since these topologies are
Fr\'{e}chet-topologies, this is equivalent to assume the operator $A$ be
linear and continuous with respect to the weak topologies $\sigma (\mathcal{S%
}_{n})$ and $\sigma (\mathcal{S}_{m})$). Let $\delta $ be the Dirac family
in $\mathcal{S}_{m}^{\prime }$ and consider the family 
\[
A^{\vee }:=(\delta _{p}\circ A)_{p\in \Bbb{R}^{m}}. 
\]
The family $A^{\vee }$ is a family in $\mathcal{S}_{n}^{\prime }$, since
each functional $A_{p}^{\vee }$ is the composition of two linear and
continuous mappings. Moreover, the family $A^{\vee }$ is of class $\mathcal{S%
}$, in fact, for every test function $\phi $ in $\mathcal{S}_{n}$ and for
every index $p$ in $\Bbb{R}^{m}$, we have 
\begin{eqnarray*}
A^{\vee }(\phi )(p) &=&A_{p}^{\vee }(\phi )= \\
&=&(\delta _{p}\circ A)(\phi )= \\
&=&\delta _{p}(A(\phi ))= \\
&=&A(\phi )(p),
\end{eqnarray*}
so that the image of the test function by the family $A^{\vee }$ is nothing
but the image of the test function under the operator $A$, i.e. 
\[
A^{\vee }(\phi )=A(\phi ), 
\]
and this image belongs to the space $\mathcal{S}_{m}$ by the choice of the
operator $A$ itself. $\blacksquare $

\bigskip

So we can give the following definition.

\bigskip

\textbf{Definition (of }$^{\mathcal{S}}$\textbf{family generated by a linear
and continuous operator). }\emph{Let }$A:\mathcal{S}_{n}\rightarrow \mathcal{%
S}_{m}$\emph{\ be a linear and continuous operator with respect to the
natural topologies of }$\mathcal{S}_{n}$\emph{\ and }$\mathcal{S}_{m}$\emph{%
\ \ (or equivalently, continuous with respect to the weak topologies }$%
\sigma (\mathcal{S}_{n})$\emph{\ and }$\sigma (\mathcal{S}_{m})$\emph{) and
let }$\delta $\emph{\ be the Dirac family in }$\mathcal{S}_{m}^{\prime }$%
\emph{. The family} 
\[
A^{\vee }:=(\delta _{p}\circ A)_{p\in \Bbb{R}^{m}} 
\]
\emph{is called \textbf{the Schwartz family generated by the operator} }$A$%
\emph{.}

\bigskip

\textbf{Remark.} We have so constructed the mapping 
\[
\left( \cdot \right) ^{\vee }:\mathcal{L}\left( \mathcal{S}_{n},\mathcal{S}%
_{m}\right) \rightarrow \mathcal{S}(\Bbb{R}^{m},\mathcal{S}_{n}^{\prime
}):A\mapsto \left( \delta _{x}\circ A\right) _{x\in \Bbb{R}^{m}}, 
\]
which we shall call \emph{the canonical representation of the operator space}
$\mathcal{L}\left( \mathcal{S}_{n},\mathcal{S}_{m}\right) $ \emph{into the
family space} $\mathcal{S}(\Bbb{R}^{m},\mathcal{S}_{n}^{\prime })$. It is
quite simple to prove that this mapping is a linear injection.

\bigskip

We shall see, as we already said, that every Schwartz family is generated by
a linear and continuous operator as specified above, or, in other terms,
that the canonical representation $\left( \cdot \right) ^{\vee }$ is a
linear isomorphism.

\bigskip

\textbf{Matrix interpretation.} As in the finite dimensional case, recalling
our interpretation of the families of distributions as generalized matrices,
we have so associated with every linear and continuous operator a
generalized matrix which is of Schwartz class.

\bigskip

\section{\textbf{The operator generated by an}$\mathcal{\ }^{\mathcal{S}}$%
\textbf{family}}

\bigskip

\textbf{Definition}\emph{\ }\textbf{(operator generated by an}$\mathcal{\ }^{%
\mathcal{S}}$\textbf{family).}\ \emph{Let }$v$\emph{\ be a family of class }$%
\mathcal{S}$\emph{\ belonging to the space }$\mathcal{S}(\Bbb{R}^{m},%
\mathcal{S}_{n}^{\prime })$\emph{. We call \textbf{operator generated by the
family} }$v$\emph{\ (or\ \textbf{associated with the family} }$v$\emph{) the
operator} 
\[
\widehat{v}:\mathcal{S}_{n}\rightarrow \mathcal{S}_{m}:\phi \mapsto v(\phi
), 
\]
\emph{sending every test function }$\phi $\emph{\ of }$\mathcal{S}_{n}$ 
\emph{into its image }$v(\phi )$\emph{\ under the family }$v$\emph{.}

\bigskip

\textbf{Example (on the Dirac family).}\emph{\ The operator (on }$\mathcal{S}%
_{n}$\emph{) generated by the Dirac family, i.e., by the ordered family }$%
\delta =(\delta _{y})_{y\in \Bbb{R}^{n}}$\emph{\ of Dirac distributions is
the identity operator of the space }$\mathcal{S}_{n}$\emph{.}

\emph{\bigskip }

\emph{Proof.} In fact, for each $y\in \Bbb{R}^{n}$, we have 
\begin{eqnarray*}
\widehat{\delta }(\phi )(y) &=&\delta _{y}(\phi )= \\
&=&\phi (y)= \\
&=&\Bbb{I}_{\mathcal{S}_{n}}(\phi )(y),
\end{eqnarray*}
for any test function $\phi $ in $\mathcal{S}_{n}$. $\blacksquare $

\bigskip

We recall that the set $\mathcal{S}(\Bbb{R}^{m},\mathcal{S}_{n}^{\prime })$,
of $^{\mathcal{S}}$families indexed by $\Bbb{R}^{m}$, is a subspace of the
vector space $(\mathcal{S}_{n}^{\prime })^{\Bbb{R}^{m}}$, of all families in 
$\mathcal{S}_{n}^{\prime }$ indexed by $\Bbb{R}^{m}$.

\bigskip

Moreover, we immediately can prove the following obvious proposition.

\bigskip

\textbf{Proposition.} \emph{For each family }$v\in \mathcal{S}(\Bbb{R}^{m},%
\mathcal{S}_{n}^{\prime })$\emph{, the operator }$\widehat{v}$\emph{\
associated with the family }$v$ \emph{is linear and the map} 
\[
\mathcal{S}(\Bbb{R}^{m},\mathcal{S}_{n}^{\prime })\rightarrow \mathrm{Hom}%
\left( \mathcal{S}_{n},\mathcal{S}_{m}\right) :v\mapsto \widehat{v} 
\]
\emph{is an injective linear operator.}

\bigskip

\textbf{Example (on the family generated by an operator). }The operator
associated with the family $A^{\vee }$ generated by a linear and continuous
operator $A$ in $\mathcal{L}\left( \mathcal{S}_{n},\mathcal{S}_{m}\right) $
is the operator $A$ itself, as can be immediately proved. In other terms we
can write 
\[
(A^{\vee })^{\wedge }=A. 
\]

\bigskip

\section{\textbf{Summability of}$\mathcal{\ }^{\mathcal{S}}$\textbf{families}
}

\bigskip

Here we give one of the fundamental definitions of Superposition Theory: the
definition of summability of a Schwartz family with respect to linear
functionals.

\bigskip

\textbf{Definition (summability of }$^{\mathcal{S}}$\textbf{families}).\emph{%
\ Let }$v$\emph{\ be a Schwartz family of tempered distributions\ belonging
to the space }$\mathcal{S}(\Bbb{R}^{m},\mathcal{S}_{n}^{\prime })$\emph{.
The family }$v$\emph{\ is said to be \textbf{summable with respect to a
linear functional }}$a$ \emph{\textbf{on the space} }$\mathcal{S}_{m}$\emph{%
\ if the composition }$u=a\circ \widehat{v}$\emph{, i.e., the linear
functional } 
\[
u:\mathcal{S}_{n}\rightarrow \Bbb{K}:\phi \mapsto a(\widehat{v}(\phi )), 
\]
\emph{is a tempered distribution (in the space }$\mathcal{S}_{n}^{\prime }$%
\emph{).}

\bigskip

It happens that a Schwartz family is summable with respect to any tempered
distribution on its index Euclidean space.

\bigskip

\textbf{Theorem (summability of }$^{\mathcal{S}}$\textbf{families}).\emph{\
Let }$v$\emph{\ be a Schwartz family of tempered distributions\ belonging to
the space }$\mathcal{S}(\Bbb{R}^{m},\mathcal{S}_{n}^{\prime })$\emph{. Then,
the family is summable with respect to every distribution in }$\mathcal{S}%
_{m}^{\prime }$\emph{; in other terms:}

\begin{itemize}
\item  \emph{\ for every tempered distribution }$a\in \mathcal{S}%
_{m}^{\prime }$\emph{, the composition }$u=a\circ \widehat{v}$\emph{, i.e.,
the linear functional} 
\[
u:\mathcal{S}_{n}\rightarrow \Bbb{K}:\phi \mapsto a(\widehat{v}(\phi )), 
\]
\emph{is a tempered distribution in the space }$\mathcal{S}_{n}^{\prime }$%
\emph{.}
\end{itemize}

\bigskip

\emph{Proof.} Let $a\in \mathcal{S}_{m}^{\prime }$ and let $\delta $ be the
Dirac family of the space $\mathcal{S}_{m}^{\prime }$. Since the linear hull 
$\mathrm{span}(\delta )$ of the Dirac family is $\sigma (\mathcal{S}%
_{m}^{\prime })$-sequentially dense in the space $\mathcal{S}_{m}^{\prime }$
(see, for example, Boccara, page 205), there is a sequence of distributions $%
\alpha =\left( \alpha _{k}\right) _{k\in \Bbb{N}}$, in the linear hull $%
\mathrm{span}(\delta )$ of the Dirac family, converging to the distribution $%
a$ with respect to the weak* topology $\sigma (\mathcal{S}_{m}^{\prime })$;
that is we have 
\[
^{\sigma (\mathcal{S}_{m}^{\prime })}\lim_{k\rightarrow +\infty }\alpha
_{k}=a. 
\]
Now, since for any natural $k$, the distribution $\alpha _{k}$ belongs to
the linear hull $\mathrm{span}(\delta )$, there exists a finite family $%
(y_{i})_{i=1}^{h}$ of points in $\Bbb{R}^{m}$ and there is a finite family
of scalars $(\lambda _{i})_{i=1}^{h}$ in the field $\Bbb{K}$ such that 
\[
\alpha _{k}=\sum_{i=1}^{h}\lambda _{i}\delta _{y_{i}}. 
\]
Consequently, by obvious calculations, we have 
\begin{eqnarray*}
\alpha _{k}\circ \widehat{v} &=&\sum_{i=1}^{h}\lambda _{i}(\delta
_{y_{i}}\circ \widehat{v})= \\
&=&\sum_{i=1}^{h}\lambda _{i}v_{y_{i}}.
\end{eqnarray*}
Hence, for every index $k\in \Bbb{N}$, the linear functional $\alpha
_{k}\circ \widehat{v}$ belongs to the space $\mathcal{S}_{n}^{\prime }$.
Now, let $s$ be the topology of pointwise convergence on the algebraic dual $%
(\mathcal{S}_{n})^{*}$, we claim that 
\[
^{s}\lim_{k\rightarrow +\infty }(\alpha _{k}\circ \widehat{v})=a\circ 
\widehat{v}. 
\]
In fact, for every test function $\phi $ in $\mathcal{S}_{n}$, we obtain 
\begin{eqnarray*}
\lim_{k\rightarrow +\infty }(\alpha _{k}\circ \widehat{v})\left( \phi
\right) &=&\lim_{k\rightarrow +\infty }\alpha _{k}\left( \widehat{v}\left(
\phi \right) \right) = \\
&=&a\left( \widehat{v}\left( \phi \right) \right) .
\end{eqnarray*}
So we have proved that the sequence of continuous linear functionals $%
(\alpha _{k}\circ \widehat{v})_{k\in \Bbb{N}}$ is pointwise convergent to
the linear functional $a\circ \widehat{v}$. Hence, by the Banach-Steinhaus
theorem (that is applicable since $\mathcal{S}_{n}$ is a barreled space),
the linear functional $a\circ \widehat{v}$ must be continuous too, i.e. $%
a\circ \widehat{v}$ should be a tempered distribution in $\mathcal{S}%
_{n}^{\prime }$. So summability of the family $v$ holds true. $\blacksquare $

\bigskip

\section{\textbf{Characterization of summability}}

\bigskip

In the following we shall denote by $\mathcal{L}(\mathcal{S}_{n},\mathcal{S}%
_{m})$ the set of all the linear and continuous operators among the two
topological vector spaces $(\mathcal{S}_{n})$ and $(\mathcal{S}_{m})$.

\bigskip

Moreover, let consider a linear operator $A:\mathcal{S}_{n}\rightarrow 
\mathcal{S}_{m}$, we say that $A$ is \emph{(topologically) transposable} if
its algebraic transpose (adjoint) $^{*}A:\mathcal{S}_{m}^{*}\rightarrow 
\mathcal{S}_{n}^{*}$ ($X^{*}$ denotes the algebraic dual of a topological
vector space $X$), defined by 
\[
^{\ast }A(a)=a\circ A, 
\]
maps the distribution space $\mathcal{S}_{m}^{\prime }$ into the
distribution space $\mathcal{S}_{n}^{\prime }$.

\bigskip

\textbf{Theorem (characterization of summability of }$^{\mathcal{S}}$\textbf{%
families}).\emph{\ Let }$v$\emph{\ be a Schwartz family of tempered
distributions\ belonging to the space }$\mathcal{S}(\Bbb{R}^{m},\mathcal{S}%
_{n}^{\prime })$\emph{. Then, the following assertions hold and they are
equivalent:}

\begin{itemize}
\item[\emph{1)}]  \emph{\ the family }$v$\emph{\ is summable with respect to
every tempered distribution }$a\in \mathcal{S}_{m}^{\prime }$\emph{;}

\item[\emph{2)}]  \emph{\ the operator }$\widehat{v}$\emph{\ is transposable;%
}

\item[\emph{3)}]  \emph{\ the operator }$\widehat{v}$\emph{\ is weakly
continuous, i.e. continuous from }$\mathcal{S}_{n}$\emph{\ to }$\mathcal{S}%
_{m}$\emph{\ with respect to the pair of weak topologies }$(\sigma (\mathcal{%
S}_{n}),\sigma (\mathcal{S}_{m}))$\emph{;}

\item[\emph{4)}]  \emph{\ the operator }$\widehat{v}$\emph{\ is continuous
from the space }$\left( \mathcal{S}_{n}\right) $\emph{\ to the space }$%
\left( \mathcal{S}_{m}\right) $\emph{.}
\end{itemize}

\emph{\bigskip }

\emph{Proof.} Note that, after the proof of the summability property (1) in
the preceding section, we have to prove only that the last three properties
are equivalent to property (1). Property (1) is equivalent to property (2)
by definition of transposable operator. Property (2) is equivalent to
property (3) because the space of linear continuous operators $\mathcal{L}((%
\mathcal{S}_{n})_{\sigma },(\mathcal{S}_{m})_{\sigma })$ is also the space
of all the transposable linear operators from the space $\left( \mathcal{S}%
_{n}\right) $ to $\left( \mathcal{S}_{m}\right) $ (see Horv\'{a}th, chap. 3,
\S\ 12, Proposition 1, page 254). Property (3) is equivalent to property
(4). In fact, since the space $\left( \mathcal{S}_{n}\right) $ is a
Fr\'{e}chet space (and then its topology coincides with the Mackey topology $%
\tau (\mathcal{S}_{n},\mathcal{S}_{n}^{\prime })$), the space $\mathcal{L}(%
\mathcal{S}_{n},\mathcal{S}_{m})$ contains the above space $\mathcal{L}((%
\mathcal{S}_{n})_{\sigma },(\mathcal{S}_{m})_{\sigma })$, of all weakly
linear and continuous operators from $\mathcal{S}_{n}$ to $\mathcal{S}_{m}$
(i.e. with respect to the pair of topologies $(\sigma (\mathcal{S}%
_{n}),\sigma (\mathcal{S}_{m}))$, see for this result Dieudonn\'{e}
Schwartz, page 91, Corollary or Horv\'{a}th, page 258, Corollary). Moreover,
the space $\mathcal{L}(\mathcal{S}_{n},\mathcal{S}_{m})$ is contained in the
space $\mathcal{L}((\mathcal{S}_{n})_{\sigma },(\mathcal{S}_{m})_{\sigma })$%
, since every continuous linear operator among two Hausdorff locally convex
topological vector spaces is weakly continuous (see proposition 3, page 256
of J. Horv\'{a}th), so the two spaces must coincide. $\blacksquare $

\bigskip

\section{\textbf{Isomorphic realization of }$\mathcal{L}\left( \mathcal{S}%
_{n},\mathcal{S}_{m}\right) $}

\bigskip

By the characterization of summability for Schwartz families we deduce the
following important corollary, which gives the canonical isomorphic
realization of the space of linear continuous operators among Schwartz test
function spaces $\mathcal{L}(\mathcal{S}_{n},\mathcal{S}_{m})$ as a space of
Schwartz families $\mathcal{S}(\Bbb{R}^{m},\mathcal{S}_{n}^{\prime })$.

\bigskip

\textbf{Corollary (of isomorphism).}\emph{\ The vector spaces }$\mathcal{S}(%
\Bbb{R}^{m},\mathcal{S}_{n}^{\prime })$\emph{\ and }$\mathcal{L}(\mathcal{S}
_{n},\mathcal{S}_{m})$\emph{\ are isomorphic. Namely, the map }$\left( \cdot
\right) ^{\wedge }$\emph{\ from the space of family }$\mathcal{S}(\Bbb{R}%
^{m},\mathcal{S}_{n}^{\prime })$\emph{\ into the space of operators }$%
\mathcal{L}\left( \mathcal{S}_{n},\mathcal{S}_{m}\right) $\emph{,
associating with each family }$v$\emph{\ its operator }$\widehat{v}$\emph{,
is a vector space isomorphism. Moreover, the inverse of the above
isomorphism is the linear mapping} 
\[
\left( \cdot \right) ^{\vee }:\mathcal{L}\left( \mathcal{S}_{n},\mathcal{S}%
_{m}\right) \rightarrow \mathcal{S}(\Bbb{R}^{m},\mathcal{S}_{n}^{\prime }) 
\]
\emph{defined by } 
\[
A\mapsto A^{\vee }:=(\delta _{p}\circ A)_{p\in \Bbb{R}^{m}}, 
\]
\emph{i.e. the canonical representation of the operator space }$\mathcal{L}%
\left( \mathcal{S}_{n},\mathcal{S}_{m}\right) $\emph{\ into the family space 
}$\mathcal{S}(\Bbb{R}^{m},\mathcal{S}_{n}^{\prime })$\emph{, which, as a
consequence, is an isomorphism too.}

\bigskip

We can give also the following definition.

\bigskip

\textbf{Definition (canonical representation of the space }$\mathcal{S}(\Bbb{%
R}^{m},\mathcal{S}_{n}^{\prime })$\textbf{).}\emph{\ The mapping} 
\[
\left( \cdot \right) ^{\wedge }:\mathcal{S}(\Bbb{R}^{m},\mathcal{S}%
_{n}^{\prime })\rightarrow \mathcal{L}\left( \mathcal{S}_{n},\mathcal{S}%
_{m}\right) :v\mapsto \widehat{v} 
\]
\emph{is called \textbf{the canonical representation of the Schwartz family
space} }$\mathcal{S}(\Bbb{R}^{m},\mathcal{S}_{n}^{\prime })$\emph{\ \textbf{%
onto the operator space} }$\mathcal{L}\left( \mathcal{S}_{n},\mathcal{S}
_{m}\right) $\emph{.}

\bigskip

\section{\textbf{Characterization of transposability}}

\bigskip

A way to see that an operator is transposable is given by the following
characterization. It is an immediate consequence of the characterization of
the $^{\mathcal{S}}$families but we want to prove it independently.

\bigskip

\textbf{Theorem.} \emph{Let }$A:\mathcal{S}_{n}\rightarrow \mathcal{S}_{m}$%
\emph{\ be a linear operator and let }$\delta $ \emph{be the Dirac family of
the space }$\mathcal{S}_{m}^{\prime }$\emph{.\ Then, the operator }$A$\emph{%
\ \ is (topologically) transposable\ if and only if, for every point }$p\in 
\Bbb{R}^{m}$\emph{, the composition }$\delta _{p}\circ A$\emph{\ is a
tempered distribution in }$\mathcal{S}_{n}^{\prime }.$

\emph{\bigskip }

\emph{Proof.} $(\Rightarrow )$ The necessity of the condition is obvious. In
fact, we have 
\[
\delta _{p}\circ A=\;^{*}A(\delta _{p}), 
\]
and so if $A$ is topologically transposable, the composition $\delta
_{p}\circ A$ is continuous. $(\Leftarrow )$ Let $a\in \mathcal{S}
_{m}^{\prime }$ be a tempered distribution; we should prove that the
composition $a\circ A$ is continuous. Since the linear hull $\mathrm{span}%
(\delta )$ is sequentially dense in the space $\mathcal{S}_{m}^{\prime }$
(see $\left[ Bo\right] $ page 205), there is a sequence of distributions $%
\left( \alpha _{k}\right) _{k\in \Bbb{N}}$ in the hull $\mathrm{span}(\delta
)$ such that 
\[
^{\sigma (\mathcal{S}_{m}^{\prime })}\lim_{k\rightarrow +\infty }\alpha
_{k}=a. 
\]
Now, since any distribution $\alpha _{k}$ lives in the hull $\mathrm{span}%
(\delta )$ there exist a finite family $(y_{i})_{i=1}^{h}$ in $\Bbb{R}^{m}$
and a finite sequence $(\lambda _{i})_{i=1}^{h}$ in $\Bbb{K}$ such that 
\[
\alpha _{k}=\sum_{i=1}^{h}\lambda _{i}\delta _{y_{i}}, 
\]
thus we have 
\begin{eqnarray*}
\alpha _{k}\circ A &=&\sum_{i=1}^{h}(\lambda _{i}\delta _{y_{i}})\circ A= \\
&=&\sum_{i=1}^{h}\lambda _{i}(\delta _{y_{i}}\circ A);
\end{eqnarray*}
hence, for every number $k\in \Bbb{N}$, the composition $\alpha _{k}\circ A$
belongs to $\mathcal{S}_{n}^{\prime }$. Let now $s$ be the topology of
pointwise convergence in the algebraic dual $\mathcal{S}_{n}^{*}$, we have 
\[
^{s}\lim_{k\rightarrow +\infty }(\alpha _{k}\circ A)=a\circ A, 
\]
in fact 
\begin{eqnarray*}
\lim_{k\rightarrow +\infty }\left( \alpha _{k}\circ A\right) \left( \phi
\right) &=&\lim_{k\rightarrow +\infty }\alpha _{k}\left( A\left( \phi
\right) \right) = \\
&=&a\left( A\left( \phi \right) \right) ,
\end{eqnarray*}
so we have that the sequence (in $\mathcal{S}_{n}^{\prime }$) of continuous
linear form $\left( \alpha _{k}\circ A\right) _{k\in \Bbb{N}}$ converges
pointwise to the linear form $a\circ A$, then, by the Banach-Steinhaus
theorem, we conclude that the composition $a\circ A$ lives also in the space 
$\mathcal{S}_{n}^{\prime }$. $\blacksquare $

\bigskip

\textbf{David Carf\`{i}}

\emph{Faculty of Economics}

\emph{University of Messina}

\emph{davidcarfi71@yahoo.it}

\end{document}